%% file: root-new.tex
\documentclass[journal]{IEEEtran}
\usepackage{ifthen}
\newboolean{showcomments}
\setboolean{showcomments}{true}
\newcommand{\aak}[1]{  \ifthenelse{\boolean{showcomments}}
{ \textcolor{blue}{(AAK says:  #1)}} {}  }
\newcommand{\jk}[1]{  \ifthenelse{\boolean{showcomments}}
{ \textcolor{red}{(JK says:  #1)}} {}  }
\newcommand{\red}[1]{  \ifthenelse{\boolean{showcomments}}
{\textcolor{red}{#1}}{#1}}
\usepackage{mathrsfs}

\usepackage{graphics} 
\usepackage{comment}
\usepackage[%
    font={small,sf},
    labelfont=bf,
    format=hang,    
    format=plain,
    margin=0pt,
    width=\textwidth,
]{caption}
\usepackage[list=true]{subcaption}
\usepackage{epsfig} 
\usepackage{mathptmx} 
\usepackage{times} 
\usepackage{amsmath} 
\usepackage{amssymb}  
\usepackage{mathtools}
\usepackage{cleveref}
\usepackage[normalem]{ulem}
\usepackage{wrapfig}

\input{macros}

\input{macros-math}
\def\texitem#1{\par\smallskip\noindent\hangindent 25pt
               \hbox to 25pt {\hss #1 ~}\ignorespaces}

\def\st{\mbox{subject to}}
\def\sys{{\rm sys}}
\newcommand{\prob}[1]{\mathbb{P}\left[ #1 \right]}
\newcommand{\rev}{r}
\newcommand{\ra}{\rightarrow}
\newcommand{\ignore}[1]{}
\allowdisplaybreaks[4]

\crefname{proposition}{Prop.}{Prop,}

\def\limsup{\lim {\rm sup}}

\ifCLASSINFOpdf
\else
\fi
\usepackage{enumitem}

\crefname{lemma}{Lemma}{}
\crefname{figure}{Figure}{}
\crefname{theorem}{Theorem}{}
\crefname{equation}{}{}

\begin{document}
%
\title{\LARGE \bf
Efficiency Fairness Tradeoff in Battery Sharing}

%
%
%

\author{Karan N. Chadha, Ankur A. Kulkarni and Jayakrishnan Nair
\thanks{Karan N. Chadha is an undergraduate at Department of Electrical Engineering,
Indian Institute of Technology
Bombay, Mumbai, India, 400076. 
{\tt\small karanchadhaiitb@gmail.com}}
\thanks{Ankur A. Kulkarni is a faculty at Systems and Control Engineering, Indian Institute of Technology
Bombay, Mumbai, India, 400076
{\tt\small kulkarni.ankur@iitb.ac.in}}%
\thanks{Jayakrishnan Nair is a faculty at Department of Electrical Engineering,
Indian Institute of Technology
Bombay, Mumbai, India, 400076
{\tt\small jayakrishnan.nair@ee.iitb.ac.in}. The authors acknowledge support from the DST grant
DST/CERI/MI/SG/2017/077.}%
}
\maketitle

\begin{abstract}
The increasing presence of decentralized renewable generation in the
power grid has motivated consumers to install batteries to save excess
energy for future use. The high price of energy storage calls for a
shared storage system, but careful battery management is required so
that the battery is operated in a manner that is fair to all and as
efficiently as possible. In this paper, we study the tradeoffs between
efficiency and fairness in operating a shared battery. We develop a
framework based on constrained Markov decision processes to study both
regimes, namely, optimizing efficiency under a hard fairness
constraint and optimizing fairness under hard efficiency
constraint. Our results show that there are fundamental limits to
efficiency under fairness and vice-versa, and, in general, the two
cannot be achieved simultaneously. We characterize these fundamental
limits via absolute bounds on these quantities, and via the notion of
\textit{price of fairness} that we introduce in this paper.
\end{abstract}

\begin{IEEEkeywords}
Smart grids, Battery management systems, Battery sharing
\end{IEEEkeywords}

\input{intro}
\input{lit-survey}

\input{model}

\input{hard_fairness}
\input{pof}
\input{soft_fairness}
\input{conclusion}
\appendices
\input{appendix}

%
\IEEEpeerreviewmaketitle

\ifCLASSOPTIONcaptionsoff
  \newpage
\fi



%
\bibliographystyle{IEEEtran}
\bibliography{ref}

%

\begin{IEEEbiography}[{\includegraphics[width=1in]{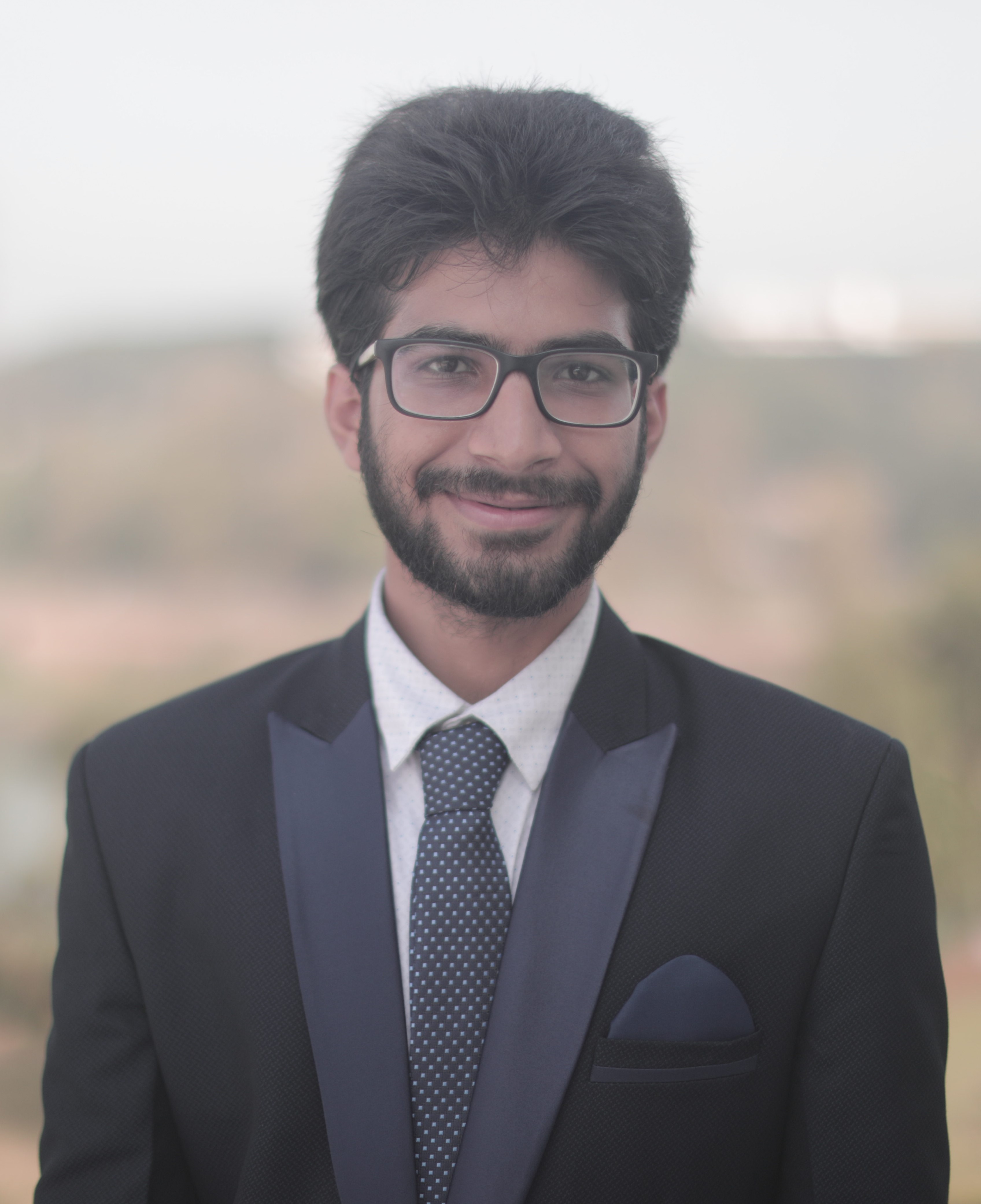}}]{Karan N. Chadha}
Karan is a dual degree student of Electrical Engineering at Indian Institute of Technology Bombay (IITB). His research interests include applied probabilty, game theory, optimization and learning theory.
\end{IEEEbiography}

\begin{IEEEbiography}[{\includegraphics[width=1in]{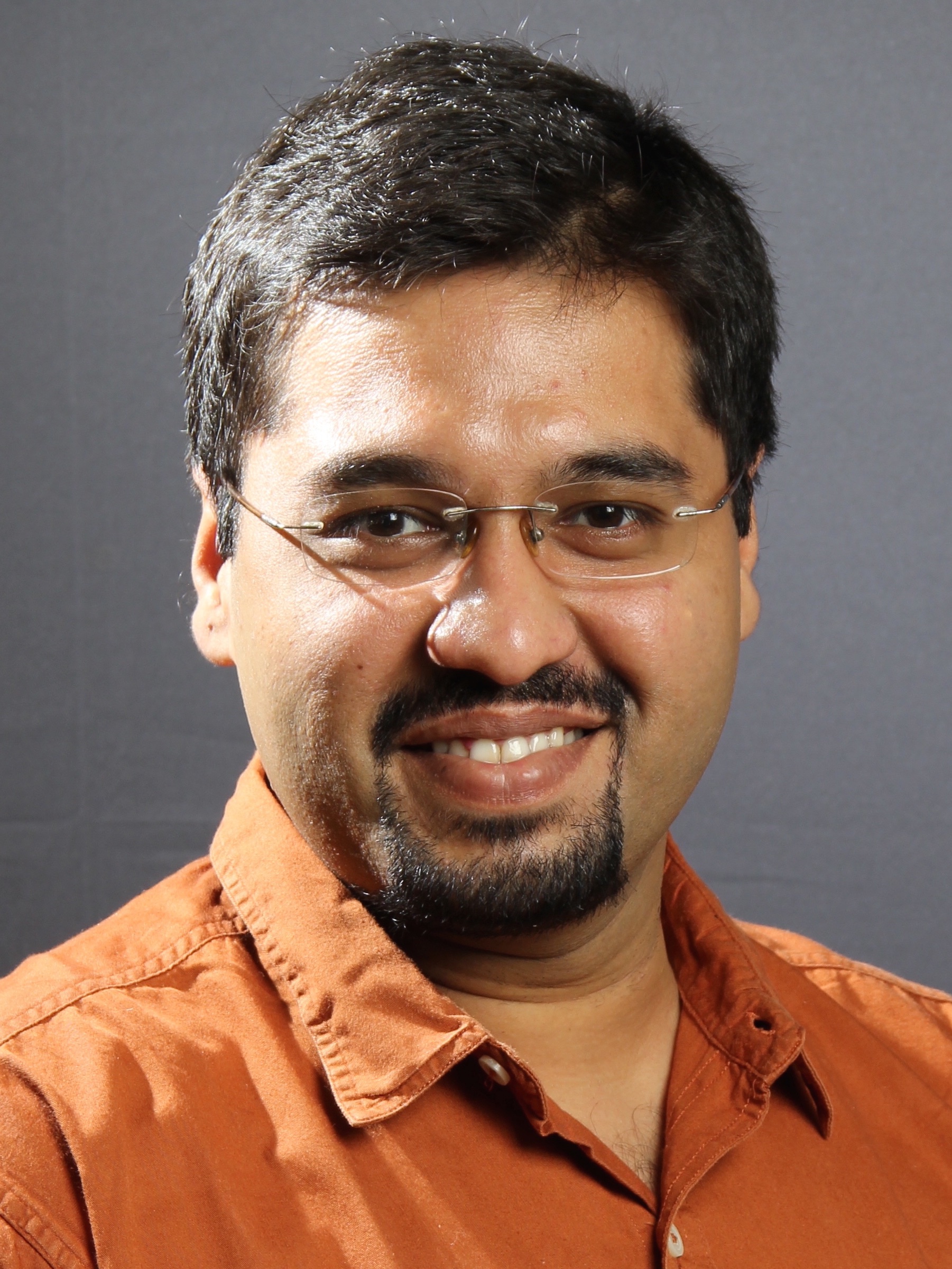}}]{Ankur A. Kulkarni}
Ankur is an Associate Professor with the Systems and Control Engineering group at Indian Institute of Technology Bombay (IITB). He received his B.Tech. from IITB in 2006, M.S. in 2008 and Ph.D. in 2010, both from the University of Illinois at Urbana-Champaign (UIUC). From 2010-2012 he was a post-doctoral researcher at the Coordinated Science Laboratory at UIUC. His research interests include information theory, stochastic control, game theory, combinatorial coding theory problems, optimization and variational inequalities, and operations research. He was an Associate (from 2015--2018) of the Indian Academy of Sciences, Bangalore, a recipient of the INSPIRE Faculty Award of the Department of Science and Technology, Government of India, 2013, Best paper awards at the National Conference on Communications, 2017, Indian Control Conference, 2018 and International Conference on Signal Processing and Communications (SPCOM) 2018, Excellence in Teaching Award 2018 at IITB  and the William A. Chittenden Award, 2008 at UIUC. He is a consultant to the Securities and Exchange Board of India on regulation of high frequency trading.
\end{IEEEbiography}

\begin{IEEEbiography}[{\includegraphics[width=1in,height=1.25in,clip,keepaspectratio]{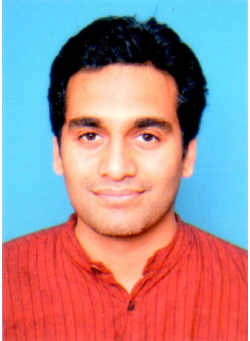}}]{Jayakrishnan
    Nair} Jayakrishnan received his BTech and MTech in Electrical Engg. (EE) from
    IIT Bombay (2007) and Ph.D. in EE from California Inst. of
    Tech. (2012). He has held post-doctoral positions at California
    Inst. of Tech. and Centrum Wiskunde \& Informatica. He is
    currently an Assistant Professor in EE at IIT Bombay. His research
    focuses on modeling, performance evaluation, and design issues in
    queueing systems and communication networks.
\end{IEEEbiography}

\end{document}

%% file: macros.tex
\usepackage{amsmath, amssymb, xspace,url}
\usepackage{epsfig}
\usepackage{longtable}
\usepackage[usenames, dvipsnames]{color}

\def\Zbb{{\mathbb{Z}}}

\def\Nbb{{\mathbb{N}}}

\def\Ebb{\mathbb{E}}

\def\Sbf{{\bf S}}

\def\Abf{{\bf A}}

\def\bfone{{\bf 1}}

\def\t{^\top}

\def\bkE{{\rm I\kern-.17em E}}
\def\bk1{{\rm 1\kern-.17em l}}
\def\bkD{{\rm I\kern-.17em D}}
\def\bkR{{\rm I\kern-.17em R}}
\def\bkP{{\rm I\kern-.17em P}}

\def\bkZ{{\bf{Z}}}
\def\bkE{{\rm I\kern-.17em E}}
\def\bk1{{\rm 1\kern-.17em l}}
\def\bkD{{\rm I\kern-.17em D}}
\def\bkR{{\rm I\kern-.17em R}}
\def\bkP{{\rm I\kern-.17em P}}
		\def\bkE{{\rm I\kern-.17em E}}
		\def\bk1{{\rm 1\kern-.17em l}}
		\def\bkD{{\rm I\kern-.17em D}}
		\def\bkR{{\rm I\kern-.17em R}}
		\def\bkP{{\rm I\kern-.17em P}}
		\def\bkY{{\bf \kern-.17em Y}}
		\def\bkZ{{\bf \kern-.17em Z}}
		\def\bkC{{\bf  \kern-.17em C}}

\let\forallnew\forall
\renewcommand{\forall}{\forallnew\ }
\def\LLR{{\rm LLR}}
\def\FC{{\rm FC}}
\def\PoF{{\rm PoF}}

\let\forall\forallnew

\def\ds{\displaystyle}

\def\b12{(\beta_1,\beta_2)}

\newcommand{\I}[1]{\mathbb{I}\{#1\}}

\def\subject{\hbox{\rm subject to}}

\def\hbar{\skew{4.2}\bar h}

\def\superstar{^{\raise 0.5pt\hbox{$\nthinsp *$}}}
\def\SUPERSTAR{^{\raise 0.5pt\hbox{$*$}}}

\def\lamstarT {\lambda^{\raise 0.5pt\hbox{$\nthinsp *$}T}}

\newlength{\noteWidth}
\long\def\notes#1{\ifinner
{\tiny #1}
\else
\marginpar{\parbox[t]{\noteWidth}{\raggedright\tiny #1}}
\fi\typeout{#1}}
 \def\notes#1{\typeout{read notes: #1}} 

\newcommand{\ie}{i.e.\@\xspace} 
\newcommand{\eg}{e.g.\@\xspace} 

\newcommand{\minimize}[1]{\displaystyle\minim_{#1}}

\newcommand{\minim}{\mathop{\hbox{\rm minimize}}}

\newcommand{\maximize}[1]{\displaystyle\maxim_{#1}}
\newcommand{\maxim}{\mathop{\hbox{\rm maximize}}}

\def\spose#1{\hbox to 0pt{#1\hss}}

\def\text #1{\hbox{\quad#1\quad}}

\makeatletter
\newcommand{\pushright}[1]{\ifmeasuring@#1\else\omit\hfill$\displaystyle#1$\fi\ignorespaces}
\newcommand{\pushleft}[1]{\ifmeasuring@#1\else\omit$\displaystyle#1$\hfill\fi\ignorespaces}
\makeatother

\def\nthinsp{\mskip -2   mu}

%% file: macros-math.tex
\newtheorem{theorem}{Theorem}[section]

\newtheorem{lemma}{Lemma}[section]

\newtheorem{proposition}[theorem]{Proposition}

\newcounter{example}
\renewcommand{\theexample}{\thesection.\arabic{example}}

\newcounter{remark}
\renewcommand{\theremark}{\thesection.\arabic{remark}}

{\begin{list}{}%
         {\setlength{\leftmargin}{#1}}%
         \item[]%
}
{\end{list}}

		\def\bsp{\begin{split}}
		\def\beq{\begin{eqnarray}}
		\def\bal{\begin{align*}}
		\def\bc{\begin{center}}
		\def\be{\begin{enumerate}}
		\def\bi{\begin{itemize}}
		\def\bs{\begin{small}}
		\def\bS{\begin{slide}}
		\def\ec{\end{center}}
		\def\ee{\end{enumerate}}
		\def\ei{\end{itemize}}
		\def\es{\end{small}}
		\def\eS{\end{slide}}
		\def\eeq{\end{eqnarray}}
		\def\eal{\end{align*}}
		\def\esp{\end{split}}
		\def\qed{ \vrule height7.5pt width7.5pt depth0pt}  

	\def\maxproblemsmall#1#2#3#4{\fbox
		 {\begin{tabular*}{0.47\textwidth}
			{@{}l@{\extracolsep{\fill}}l@{\extracolsep{6pt}}l@{\extracolsep{\fill}}c@{}}
				#1 & $\maximize{#2}$ & $#3$ & $ $ \\[5pt]
					 & $\subject\ $    & $#4$ & $ $
			\end{tabular*}}
			}

		\def\problemsmall#1#2#3#4{\fbox
		 {\begin{tabular*}{0.48\textwidth}
			{@{}l@{\extracolsep{\fill}}l@{\extracolsep{6pt}}l@{\extracolsep{\fill}}c@{}}
				#1 & $\minimize{#2}$ & $#3$ & $ $ \\[5pt]
					 & $\subject\ $    & $#4$ & $ $
			\end{tabular*}}
			}

	\def\cp2problem#1#2#3#4{\fbox
		 {\begin{tabular*}{0.9\textwidth}
			{@{}l@{\extracolsep{\fill}}l@{\extracolsep{6pt}}l@{\extracolsep{\fill}}c@{}}
				#1 & & $#4 $ 
			\end{tabular*}}}

		\def\bkE{{\rm I\kern-.17em E}}
		\def\bk1{{\rm 1\kern-.17em l}}
		\def\bkD{{\rm I\kern-.17em D}}
		\def\bkR{{\rm I\kern-.17em R}}
		\def\bkP{{\rm I\kern-.17em P}}
		
		\def\bkZ{{\bf{Z}}}

\newcommand {\beeq}[1]{\begin{equation}\label{#1}}
\newcommand {\eeeq}{\end{equation}}
\newcommand {\bea}{\begin{eqnarray}}
\newcommand {\eea}{\end{eqnarray}}

\def\texitem#1{\par\smallskip\noindent\hangindent 25pt
               \hbox to 25pt {\hss #1 ~}\ignorespaces}



%% file: intro.tex
\section{Introduction}
The falling cost of solar panels and incentives for renewable
generation has created hope for massively decentralized power
generation. Islanded microgrids can install wind generators,
individual households, housing societies, and industries can generate
their own energy via solar panels, and reduce both, their reliance on
the grid and their carbon footprint. However, unpredictability and
unreliability of renewable energy also necessitate battery backup to
store excess energy, so that it can be utilized at times when this
generation is low. With adequate battery backup, one can envision a
future where there is widespread deployment of renewable generation,
leading to a higher penetration of clean energy sources. A key
component in realizing this dream is the cost of the battery. While
solar panels have become inexpensive in the recent past, batteries
prices continue to be stiff and significant.

The high cost of battery storage has motivated communities to look for
sharing contracts that provide users access to a common battery,
leading to sharing of the overall costs. The shared battery is
expected to serve as collective storage accessible by multiple users
that can be utilized to store excess generation, which users may want
to withdraw at times of deficit.
However, like anyone who has used a shared refrigerator at a dorm
would realize, sharing of a storage space brings with it concerns
about the inequity of access and the possibility of other users
`stealing' one's content. Clearly, gate-keeping is required so that
users access the battery in a manner that is fair to all, and yet the
battery is operated in an efficient manner that fully exploits the
investment made in installing it. We are motivated by this question of
fair and efficient operation of a shared battery.

For the success of the sharing arrangement, it should ideally allow
each user to access the common battery as if it is wholly owned by the
user without jeopardizing the access of other users. But if users have
disparate generation characteristics, excessive charging by one user
may leave another user with less opportunity to deposit his energy,
diminishing equity of access. At the other extreme, it is evident that
without a formal algorithm to do the gate-keeping, some users could
draw far more energy than they injected, thereby cannibalizing the
energy of other users.

Moreover, it seems plausible that gate-keeping may compel the battery
operator to decline requests for injection even while there is room in
the battery, or for withdrawal of energy even there is energy in the
battery. As a consequence, maintaining \textit{fairness} may come in
the way of \textit{efficient} battery operation. In this paper, we
develop a framework to quantify the tradeoffs between efficiency and
fairness in the operation of a shared battery. Our work shows that
there is a precise sense in which the requirement of fairness
fundamentally compromises efficiency, and in general, the two cannot
be achieved simultaneously. We also quantify the tradeoff by
characterizing the maximum efficiency achievable under fairness, and
the maximum fairness achievable under efficiency.

We consider a setting where $N$ users with individual stochastic net
generation access a shared battery. The net generation may be positive
(in which case the user wishes to inject energy) or negative (draw
energy). We posit that a battery management algorithm decides the
amount of energy to accept from injecting users and to provide to
demanding users.  We adopt a general principle of fairness that no
user should draw more energy from the battery than it has injected
into it.  However, since energy units added by users are fungible, in
the sense that they are indistinguishable from those added by other
users, the above principle need not be applied sequentially. The
battery management algorithm can allow users an `overdraft', where
users can for some intervals, draw more energy from the battery than
they have injected thus far, so long as the energy is `returned' at a
later stage.  Specifically, in our notion of fairness, we ask that for
each user, the long-run time average of the energy drawn from the
battery by the user be no more than the long-run time average of the
energy injected. Our notion of efficiency is the \textit{loss of load
  rate} (LLR). The LLR of user $i$ is its long run time average of
unmet demand.

We first study the problem of maximizing the efficiency subject to
hard constraints on fairness and derive fundamental limits on the
minimum total LLR attainable under fairness constraints. We show that
when there is at least one user that is \textit{net demanding}, \ie,
has negative steady-state net generation on average; the total LLR
remains bounded away from zero for any battery size. Remarkably, this
is true even when the total steady-state average net generation of all
users is positive. In contrast, in the absence of the fairness
constraint, the total LLR, in this case, would decrease exponentially
with battery size. When all users are net generative (\ie, with
positive steady state net generation on average), the LLR even with
the fairness constraint decreases to zero exponentially with the
battery size.

Given that the fairness constraint limits the efficient use of the
battery, we then ask the \textit{price of fairness}. This is defined
as the ratio of the minimum attainable total LLR with the fairness
constraint, to that attained without the fairness constraint. Our
fundamental limits already show that this ratio approaches infinity
for large battery size when at least one user is net demanding, but
the system is net generative on the whole. However, remarkably, we
show that the price of fairness can be arbitrarily large even all
users are net generative.


Finally, we study the maximum fairness attainable under efficient
battery use. We find numerically that efficient battery operation
often compromises fairness for at least one user. Moreover, that
departure from the utopian fairness persists even when one increases
the battery size. This shows that the conflict between efficiency and
fairness is, in a sense, fundamental and that a larger battery does
not remedy it. These results provide a natural motivation for devising a market for 
energy wherein a user can inject its excess energy into the battery, which can then be 
sold to other users. The precise formulation for such a market is a point for future study.

%% file: lit-survey.tex
\subsection*{Literature Survey}

There is recent literature on the effects of energy/battery sharing in
power systems. The effects of sharing energy are studied from a
(non-cooperative) game-theoretic perspective in
\cite{kalathil2019sharingelec}. The strategic decision is to choose
the investment in individual battery storage systems. On the other
hand, \cite{chakraborty2018coalationalsharing} views this problem from
a coalition game perspective.
Investment in a shared battery and the corresponding allocation scheme
is studied in \cite{wu2016communitystorage}. All these models assume
the net generation observed by the battery to be independent and
identically distributed random variables (i.i.d.) across time. This
seriously limits the applicability of the system as the net generation is
generally dependent across time. In contrast, we consider Markov
models for net generation, which can better model the behaviour of
renewable energy sources.

Markov energy generation models for scheduling of energy storage have
been considered in \cite{zhou2018windstoragemdp},
\cite{kim2011optenergy}. \cite{zhou2018windstoragemdp} models the
system as a finite-horizon Markov decision process (MDP) where the
decision of buying or selling electricity is made by an operator. They
provide threshold based heuristic policies and quantify their
suboptimality. \cite{korpaas2003opersizingstorage} studies the problem
of scheduling and operation of energy storage for wind power
plants. \cite{vandevan2012optcontenergy} considers an infinite horizon
discounted MDP based model to determine the optimal amount of battery
to be bought from the grid to minimize the cost given a battery of
fixed capacity. These papers focus on providing structural results of
optimal policies and/or heuristic algorithms. On the other hand, we
consider a model in which we minimize the expected loss of load if
all demand is to be satisfied by only the renewable generation and the
common battery. Our main distinction is the introduction of a notion
of fairness, which, to the best of our knowledge, has not been done
before in the scheduling of energy storage. Using this notion, we
study the tradeoff between efficiency and fairness in the operation of
shared energy storage systems.

%% file: model.tex
\section{Model and Preliminaries}

\subsection{Notation}
In general, we denote random variables by capital letters and the
value taken by a random variable is denoted by the corresponding lower
case letters. $\prod_{i \in I}V_i$ denotes the Cartesian product of
sets $V_i, \ i \in I$. $\{u,\ldots,v\}$ is the set of integers from
$u$ to $v$. $[n]$ denotes the set $\{1,\hdots,n\}$.

\subsection{Model}

Consider $n$ users $\{1,\hdots, n\}$ equipped with stochastic net
generation evolving in discrete time. The net energy generation of
user $i$ at time $t$ be denoted by $X_i(t) \in S_i.$ A positive value
of $X_i(t)$ indicates a net surplus at time~$t$ (i.e., user~$i$
generated more energy than she consumed) whereas a negative value of
$X_i(t)$ indicates a net deficit at time $t$ (i.e., user~$i$ demanded
more energy than she generated). We assume an arbitrary positive
granularity with which energy generation/demand is measured; this
granularity is taken to be unity without loss of generality, so that
$S_i \subset \Zbb.$
Let $s^+_{i}$ denote the maximum energy user~$i$ injects and
$-s^-_{i}$ denote the maximum energy user~$i$ demands, i.e., $s^+_{i}
= \max_{s \in S_i} s$ and $s^-_{i} = \min_{s \in S_i} s$. To avoid
degenerate scenarios, we assume that $s^+_{i} > 0$ and $s^-_{i} < 0.$

Let $X(t) =\{X_1(t),\hdots,X_n(t)\}$. We assume that $X(\cdot)$ is an
irreducible discrete time Markov chain (DTMC) over a finite state
space $S \subseteq \prod_{i \in [n]} S_i.$
Note that we do not assume here that the net generation processes
associated with the individual users are independent.\footnote{It is
  straightforward to generalize our results to the more general
  setting where the vector of net generations is itself a function of
  an abstract background Markov process. This generalization allows
  for an arbitrary state space desciption that might, for example,
  incorporate history and/or weather information.} We make the
assumption that $\prob{X(t+1)=s\ |\ X(t)=s} > 0$ for all $s \in S.$
Note that this ensures aperiodicity of $X(\cdot).$
%

Let $\pi = (\pi(s),s \in S)$ denote the stationary distribution of the
DTMC $X(\cdot).$ We define the \emph{drift} $\Delta_i$ associated with
user~$i$ as her steady state average net generation: $$\Delta_i =
\sum_{s \in S} s_i \pi(s),$$ 
where $s=(s_1,\hdots,s_n).$ 
User~$i$ is said to be \emph{net
  generative} if $\Delta_i > 0$ and \emph{net demanding} if
$\Delta_i<0.$ The system drift $\Delta$ is defined the sum of the user
drifts, i.e., $\Delta = \sum_{i \in [n]} \Delta_i.$

A common battery with capacity $b_{\max} \in \Nbb$ is shared between
the users. The battery occupancy at time $t$ is a random variable
denoted by $B(t)$. $(X(t),B(t))$ is a controlled Markov process
evolving over
\begin{equation}\label{eqn:state-space}
\Sbf = S \times [b_{\max}].
\end{equation}
Battery dynamics are given by
\begin{equation}\label{eqn:bat-dyn}
B(t+1) = B(t) + \sum_{i\in [n]} A_i(t),
\end{equation}
where $A_i(t)$ denotes the energy accepted from user $i$ (when $A_i(t)
\geq 0$) or the energy supplied to user $i$ (when $A_i(t) < 0$).  We
assume that the actions $A_i(t), i\in [n]$ at each time instant~$t$
are chosen by the battery operator as a function of the state history
$(X(s),B(s)),$ $s \leq t$.

The battery management algorithm is constrained as follows. In state $(x,b) \in \Sbf$, the space of allowable actions, denoted by $A(x,b),$
is given by: 
\begin{align}\label{eqn:act-space}
\nonumber A(x,b) &= \biggl\{ (a_1,\ldots,a_n)\ |\  a_i \in  \begin{cases}
  \{0,\ldots,x_i\}, & {\rm if} \ x_i \geq 0  \\
  \{x_i,\ldots,0\},  & {\rm if} \ x_i < 0
\end{cases}, \\
& \quad \quad 0 \leq b + \sum_{i\in [n]} a_i \leq b_{\max} 
\biggr\} \\
\Abf &= \bigcup_{(x,b) \in \Sbf} A(x,b) \nonumber
\end{align}
The above constraints restrict the amount of energy supplied to a user
by the amount demanded and similarly restricts the amount of energy
accepted from a user by her net surplus. In addition, the actions must
also respect the capacity constraints of the battery, and that the
battery cannot be discharged below~0.

\subsection{Loss of Load Rate (LLR)}
Our notion of efficiency is defined by the loss of load rate.  The
loss of load rate ($\LLR$) for user $i,$ denoted by $\LLR_i,$ is
defined as:
\begin{equation}\label{eqn:llr_i}
\LLR_i \coloneqq \lim_{T\rightarrow \infty} \frac{1}{T} \sum_{t=0}^T \Ebb[ \I{X_i(t)<0}(A_i(t)-X_i(t))]
\end{equation}
Note that $\LLR_i \geq 0,$ since $A_i(t) \geq X_i(t)$ when $X_i(t)< 0$
from~\cref{eqn:act-space}. $\LLR_i$ captures the long run average rate
of unmet demand for user~$i.$
For a system of $n$ users, we define the $\LLR$ of the system to be
the sum of $\LLR$s of the individual users, \ie, $\LLR_\sys =
\sum_{i \in [n]} \LLR_i$.

%% file: hard_fairness.tex
\section{Maximizing Efficiency with Hard Constraints on Fairness}
\label{sec:hard_fairness}

In this section, we study the tradeoff between efficiency and fairness
by putting a hard constraint on fairness and optimizing efficiency
under this constraint. The notion of fairness we consider is that for
each user, the time-averaged amount of energy drawn from the battery
is at most the time-averaged amount of energy injected into the
battery by that user. Efficiency is measured using the LLR. We
formulate this problem as a constrained Markov decision process, and
based on this formulation, derive fundamental limits on the efficiency
achievable under fairness constraints.

\subsection{Constrainted Markov decision process formulation}
In this section, we impose a set of service constraints that the
battery management algorithm must satisfy, in order to impart fairness
in battery scheduling.
\begin{equation}\label{eqn:fair-cons}
\FC_i: \qquad \lim_{T \rightarrow \infty} \frac{1}{T}\sum_{t=0}^T \Ebb[A_i(t)] {\geq} 0.
\end{equation}
$\FC_i$ asks that the amount of energy we provide to each user is at
most the amount of energy the user injects into the battery.  We call
$\lim_{T \rightarrow \infty} \frac{1}{T}\sum_{t=0}^T \Ebb[A_i(t)] $
the {\it net contribution} (C$_i$) of source $i$.  We impose that the
battery management algorithm satisfy $\FC_i$ for all $i \in [n]$,
collectively referred to as \textit{fairness constraints}.

We consider the problem of operating the system under fairness
constraints so that $\LLR_\sys$ is minimized. This problem is an
instance of a constrained Markov decision process (we refer the reader
to~\cite{altman99constrainedmdp} for a survey), or CMDP with the
average cost criterion. The state space is given
by~\cref{eqn:state-space}, and the set of allowable actions for each
state is given by~\cref{eqn:act-space}. The problem, for a fixed
initial distribution on the state space, is summarized below:
\\\\ \problemsmall{(P)}{\phi}{\underset{{T\rightarrow \infty}}{\lim}
  \frac{1}{T} \underset{t=0}{\overset{T}\sum} \underset{{i\in
      [N]}}\sum \Ebb[ \I{X_i(t)<0}(A_i(t)-X_i(t))]}{\underset{T
    \rightarrow \infty}{\lim}
  \frac{1}{T}\underset{t=0}{\overset{T}\sum} \Ebb[A_i(t)] {\geq} 0
  \qquad \forall i \in [n].}

Here the decision variable is $\phi$, a randomized, history-dependent
policy. A randomized policy is a sequence of functions $\phi(t), t\in
\Nbb$ that map the set of histories until time~$t$ to probability
distributions on the set of available actions at time~$t$.  The set of
available actions in state $(x,b) \in \Sbf$ is given by
\eqref{eqn:act-space}. A deterministic policy maps the set of
histories to a specific available action at each time $t\in \Nbb.$ A
policy is Markov if it depends only on current state $(X(t),b(t))$ and
it is stationary if $\phi(t)$ does not depend on $t$.

The CMDP (P) is tractable since it can be solved by solving an
equivalent linear problem (LP) (see
\cite{altman99constrainedmdp,PutermanMDP}); we do not provide the
details of this reduction here due to space constraints, though we do
present the LP reduction for the CMDP posed in
Section~\ref{soft_fairness} for optimization of fairness given a hard
constraint on efficiency. Note also that the CMDP (P) is defined for a given initial distribution 
on $\Sbf$; as such its optimal value depends on this initial distribution. However, the choice of this initial distribution
does not affect the results we derive in this paper.

\subsection{Fundamental limits}
Denote the set of net demanding and net generating sources by 
\[\mathscr{D} := \{ i \in [n] | \Delta_i< 0\}, \quad \mathscr{G} := \{ i \in [n] | \Delta_i> 0\}.  \]
The following theorem provides a lower bound on the objective of~(P)
under any policy.

\begin{theorem}\label{thm:netdemlb}
 \[
  \LLR_\sys \geq \sum_{i \in \mathscr{D}}(-\Delta_i).
 \]
\end{theorem}
\begin{IEEEproof}
  It suffices to prove the stated lower bound on~$\LLR_\sys$ under any
  feasible policy for~(P).
  \begin{align*}
    \LLR_i &= \lim_{T\rightarrow \infty} \frac{1}{T} \sum_{t=0}^T
    \Ebb[ \I{X_i(t)<0}(A_i(t)-X_i(t))] \\
    &= \lim_{T\rightarrow  \infty} \frac{1}{T} \sum_{t=0}^T \Ebb[\I{X_i(t)<0} A_i(t)] \\
    &\quad + \lim_{T\rightarrow  \infty} \frac{1}{T} \sum_{t=0}^T \Ebb[\I{X_i(t)>0} X_i(t)]- \lim_{T\rightarrow  \infty} \frac{1}{T} \sum_{t=0}^T \Ebb[X_i(t)]\\
    &\stackrel{(a)}\geq \lim_{T\rightarrow  \infty} \frac{1}{T} \sum_{t=0}^T \Ebb[\I{X_i(t)<0} A_i(t)] \\
    &\quad + \lim_{T\rightarrow  \infty} \frac{1}{T} \sum_{t=0}^T \Ebb[\I{X_i(t)>0} A_i(t)]- \lim_{T\rightarrow  \infty} \frac{1}{T} \sum_{t=0}^T \Ebb[X_i(t)] \\
    &= \lim_{T\rightarrow  \infty} \frac{1}{T} \sum_{t=0}^T \Ebb[A_i(t)] - \lim_{T\rightarrow  \infty} \frac{1}{T} \sum_{t=0}^T \Ebb[X_i(t)] \\
    &\stackrel{(b)}\geq - \lim_{T\rightarrow  \infty} \frac{1}{T} \sum_{t=0}^T \Ebb[X_i(t)] = -\Delta_i
  \end{align*}
  The inequality $(a)$ follows since $A_i(t) \leq X_i(t)$ when $X_i(t)
  > 0.$ The inequality $(b)$ is a consequence of our fairness
  constraint. The statement of the lemma now follows by summing over
  all net demanding sources.
  \ignore{
    Consider user $i$. The fairness constraint (FC$_i$) which
   says $\sum_{u = s^-_{i}}^{s^+_{i}} \rho(a_i = u)u \geq 0 $ can be
   rewritten as:
   \begin{align}
     \nonumber \sum_{s_i \in S_i: s_i < 0}\sum_{u = s_i}^{0} \rho(x_i = s_i,a_i = u)(u-s_i) + \\ 
     \label{eqn:pf-low1}\sum_{s_i \in S_i: s_i \geq 0}\sum_{u = 0}^{s_i} \rho(x_i = s,a_i = u)(u-s_i) 
     & \geq -\sum_{s_i \in S_i}s_i\rho(x_i = s_i) 
   \end{align}
   \cref{eqn:pf-low1} follows by adding $-\sum_{s_i \in S}s\rho(x_i =
   s)$ on both sides and splitting the sum on left hand side into two
   parts depending on whether the state is injecting (positive) or
   demanding (negative). Now because the second term in the LHS of
   \cref{eqn:pf-low1} is non-positive, we get \cref{eqn:pf-low2}
   \begin{align}
     \label{eqn:pf-low2} \sum_{s_i \in S_i: s_i < 0}\sum_{u = s_i}^{0} \rho(x_i = s_i,a_i = u)(u-s_i) & \geq \sum_{s_i \in S_i}(-s_i)\rho(x_i = s_i)
   \end{align}
   The RHS of \eqref{eqn:pf-low2} is $-\Delta_i$.  Thus, inspecting
   the objective of (LP) we get that $ \LLR_\sys \geq \sum_{i \in
     \mathscr{D}}(-\Delta_i), $ as required.}
\end{IEEEproof}

The above bound shows that $\LLR_\sys$ remains bounded away from zero
when $\mathscr{D} \neq \emptyset$ \textit{for any battery size}. We
next consider the case when all sources are net generative, \ie,
$\mathscr{G}=[n]$. We show that in this case the optimal $\LLR_\sys$
goes down to zero exponentially with increasing battery size. To show
this, we upper bound the optimal $\LLR_\sys$ by an expression which
exponentially decays to 0. For this, we extensively use
\cref{lem:netgensingle} from Appendix \ref{app:singlesource}. It shows
that the optimal LLR for a single net generating user decays
exponentially in battery size. Using this, we now show that the
exponential decay of the optimal LLR is true even for multiple
sources, if all of them are net generating.

\begin{theorem}
If $\mathscr{G}=[n]$, we have
 \[
  \underset{b_{\max} \rightarrow \infty}{\limsup}\frac{\log(\LLR_o)}{b_{\max}} = -c,
 \]
 where $\LLR_o$ denotes the optimal $\LLR$, $b_{\max}$ is the battery
 size and $c$ is a positive constant.
\end{theorem}
\begin{IEEEproof}
We upper bound the solution of (P) by considering a particular
policy. We divide the battery into $n$ equal chunks of size
$\frac{b_{\max}}{n},$ allocating one chunk to each user. Consider the
policy that optimally schedules each user using only the battery chunk
allocated to her; as we mention in Appendix~\ref{app:singlesource},
simple greedy operation is optimal when a battery (chunk) is used by
used by a single user. Moreover, it is easy to see that this policy is
fair.

Under the above policy, let $\LLR_i$ denote the loss of load rate of
source~$i$. Then, the $\LLR_\sys$ of the combined system under this
policy is $\sum_{i = 1}^{n}\LLR_i$. We know from
\cref{lem:netgensingle} that
 \begin{equation}\label{eqn:llr-i}
  \lim_{b_{\max} \rightarrow \infty}\frac{\log(\LLR_i)}{b_{\max}} = -c_i,  
 \end{equation}
 for some positive constant $c_i$. Without loss of generality, let
 $c_1 \leq c_2 \leq \ldots \leq c_n$. It is not hard to see now that 
 \begin{equation*}\label{eqn:llr-o}
   \underset{b_{\max} \rightarrow \infty}{\limsup}\frac{\log(\LLR_\sys)}{b_{\max}} = - c_1,   
 \end{equation*}
 which implies the statement of the theorem.
\end{IEEEproof}

%% file: pof.tex
\section{Price of Fairness}\label{sec:pof}

In this section, we analyse the efficiency implications of the
fairness constraint~\eqref{eqn:fair-cons} in the LLR optimization
formulation (P) introduced in the previous section. We introduce the
notion of Price of Fairness (PoF), which is the ratio of the optimal
LLR subject to the fairness constraint, to the optimal LLR without
this constraint. The main conclusion of this section is that the PoF
can be arbitrarily large, even when all users are net generative. This
means that imposing a strong fairness constraint in battery sharing
can result in a substantial loss of efficiency. In the following
section, we address this issue by proposing an alternative formulation
that seeks to maximize fairness subject to maximal efficiency.

To define PoF, let $\LLR_e$ denote the optimal loss of load rate
without the fairness constraint, i.e., 
\[
\LLR_e \coloneqq \min_{\phi} \lim_{T\rightarrow \infty} \frac{1}{T}
\sum_{t=0}^T \sum_{i = 1}^{n} \Ebb[ \I{X_i(t)<0}(A_i(t)-X_i(t))].
\]
We refer to policies that solve the above optimization as
\emph{efficient} policies. It is not hard to see that the optimal loss
of load rate is achieved by any greedy policy that~(i)~always charges
the battery whenever energy is available, and~(ii)~always meets user
demand whenever feasible. Formally, any policy that chooses actions
that satisfy the following conditions is efficient.
\begin{enumerate}[start=1,label={E\arabic*.}]
\item If $0 \leq b+\sum_{i \in [n]} x_i \leq b_{\max},$ then $a = x.$
\item If $b+\sum_{i \in [n]} x_i > b_{\max},$ then $a_i = x_i$ for all
  $i$ satisfying $x_i < 0,$ and $b+\sum_{i \in [n]} a_i = b_{\max}.$
\item If $b+\sum_{i \in [n]} x_i < 0,$ then $a_i = x_i$ for all $i$
  satisfying $x_i > 0,$ and $b+\sum_{i \in [n]} a_i = 0.$
\end{enumerate}
Let the space of actions $(a_1,\hdots,a_n)$ that satisfy the above efficiency conditions in state $(x,b)$ 
be denoted by $A_e(x,b)$ and the history dependent 
policies that satisfy these action constraints be denoted by $\Phi_e.$

The price of fairness (PoF) is defined as
\begin{equation}\label{eqn:pofdef}
  \PoF = \frac{\LLR_o}{\LLR_e}.
\end{equation}

To analyse the PoF, we first consider the case where there is at least
one net demanding source. In this case, if the system drift is itself
negative, it is easy to see that $\LLR_o$ is bounded away from zero
(using the same argument as in the proof of
Theorem~\ref{thm:netdemlb}). Thus, the more interesting case is
$\Delta > 0.$ In this case, the following lemma shows that the PoF
grows to infinity as $b_{\max} \ra \infty.$
\begin{lemma}
  \label{lemma:pof_nd}
  Suppose that $|\mathscr{D}| \geq 1$ (i.e., $\Delta_i < 0$ for
  some~$i$) but the system is net generative (i.e., $\Delta >
  0$). Then PoF $\ra \infty$ as $b_{\max} \ra \infty.$
\end{lemma}
\begin{IEEEproof}
  By Theorem~\ref{thm:netdemlb}, $\LLR_o$ is bounded away from zero
  for any value of $b_{\max}.$ However, $\LLR_e$ decays exponentially
  with $b_{\max};$ this is because under an efficient policy, one can
  think of the resulting $\LLR$ as arising from a single user with net
  generation process $\sum_{i \in [n]} X_i(\cdot).$ The exponential
  decay of $\LLR_e$ with $b_{\max}$ now follows from
  Theorem~\ref{lem:netgensingle} in the appendix.
\end{IEEEproof}

Next, we consider the case where all users are net generating. In this
case, one might expect that a strict fairness constraint does not
affect efficiency.
%
%
However, as the following lemma shows, the PoF can be arbitrarily
large even when all sources are net generating.
\begin{lemma}
  Given any $M > 0,$ one can construct a system instance where
  $\mathscr{G}=[n],$ such that the PoF $\geq M.$
\end{lemma}
\begin{IEEEproof}
  We will construct an instance with 2 users that satisfies PoF $\geq
  M$ for any arbitrary positive threshold $M.$ The full details of the
  construction are rather cumbersome; we provide a sketch of the
  argument below.

  For users 1 and 2, we construct two independent net generation
  processes such that the LLR decay rate corresponding to each user
  operating alone (the framework analysed in the appendix) is
  distinct, i.e., $\lambda_1 < \lambda_2.$ In this case, it can be
  shown (using the precise charactization of the decay rate from large
  deviations theory) that the decay rate $\lambda_e$ associated with
  battery sharing between the users under an efficient policy
  satisfies $\lambda_e \in (\lambda_1,\lambda_2).$ Denoting the
  standalone $LLR$ of user $i$ by $\LLR_{o,i},$ this
  implies $$\lim_{b_{\max} \ra \infty} \frac{\LLR_{o,1}}{\LLR_e} =
  \infty.$$ Thus, there exists $\bar{b}_{\max} > 0$ such
  that $$\frac{\LLR_{o,1}(\bar{b}_{\max})}{\LLR_e(\bar{b}_{\max})} >
  2M.$$ From hereon, we freeze the battery size to this value
  $\bar{b}_{\max},$ so the dependence of LLR on $\bar{b}_{\max}$ will
  be suppressed.

  Let $\delta = \LLR_{o,1}.$ We now perturb the net generation process
  of user~2 as follows: $\tilde{X}_2(s) = X_2(s) + a,$ for $a \geq 0$
  such that $$\sum_{s \in S} s \pi(s) \I{\tilde{X}(s) < 0} \in (-
  \delta/2,0).$$ Note that quantities in the perturbed system are
  represented with a \textit{tilde} accent. In the perturbed system,
  user~2 is more generating than in the original system, and the rate
  at which it demands energy from the system is at most $\delta/2.$

  It is easy to see that $\widetilde{\LLR}_e \leq \LLR_e.$ Moreover,
  in the perturbed system operated under the fairness constraint, the
  rate at which energy is accepted from user~2 is at most $\delta/2.$
  This means that the reduction in the LLR of user~1 is at most
  $\delta/2$ relative to standalone operation,
  i.e., $$\widetilde{\LLR}_o \geq \LLR_{o,1} - \frac{\delta}{2} \geq
  \frac{1}{2}\LLR_{o,1}.$$ Thus, $$\widetilde{\mathrm{PoF}} =
  \frac{\widetilde{\LLR}_o}{\widetilde{\LLR}_e} \geq \frac{1}{2}
  \frac{\LLR_{o,1}}{\LLR_e} =M.$$ This completes the proof.
\end{IEEEproof}

\newsavebox{\smlmatdemfir}
\savebox{\smlmatdemfir}{$\left[\begin{smallmatrix}0.4&0.6\\0.4&0.6\end{smallmatrix}\right]$}

\newsavebox{\smlmatgenfir}
\savebox{\smlmatgenfir}{$\left[\begin{smallmatrix}0.6&0.4\\0.6&0.4\end{smallmatrix}\right]$}

\newsavebox{\smlmatgensec}
\savebox{\smlmatgensec}{$\left[\begin{smallmatrix}0.7&0.3\\0.7&0.3\end{smallmatrix}\right]$}

\newsavebox{\smlmatgenthird}
\savebox{\smlmatgenthird}{$\left[\begin{smallmatrix}0.5&0.5\\0.6&0.4\end{smallmatrix}\right]$}

\newsavebox{\smlmatgenextreme}
\savebox{\smlmatgenextreme}{$\left[\begin{smallmatrix}0.95&0.05\\0.95&0.05\end{smallmatrix}\right]$}

\newsavebox{\smlmatgenmild}
\savebox{\smlmatgenmild}{$\left[\begin{smallmatrix}0.51&0.49\\0.51&0.49\end{smallmatrix}\right]$}

\newsavebox{\smlmatgenone}
\savebox{\smlmatgenone}{$\left[\begin{smallmatrix}0.9&0.1\\0.9&0.1\end{smallmatrix}\right]$}

\newsavebox{\smlmatgentwo}
\savebox{\smlmatgentwo}{$\left[\begin{smallmatrix}0.55&0.45\\0.55&0.45\end{smallmatrix}\right]$}

\begin{figure}
\centering
\subcaptionbox{$ P_1  = $ \usebox{\smlmatgenextreme}, \\ $ P_2  = $\usebox{\smlmatgenmild} \label{fig:pofgeneratingexp}}{\includegraphics[width=0.49\columnwidth, trim=2in .05in 2.45in .15in, clip=true]{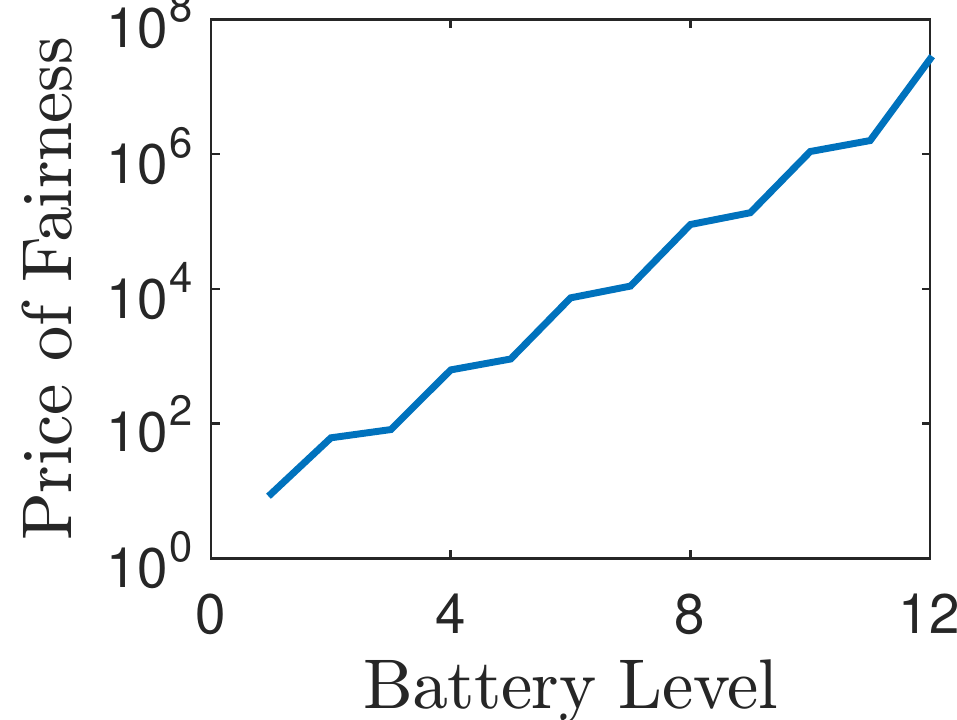}}
\subcaptionbox{$ P_1  = $ ~\usebox{\smlmatgenone}, \\ \hspace*{5mm}$ P_2  = $ ~\usebox{\smlmatgenmild},  $ P_3  = $ \usebox{\smlmatgentwo} \label{fig:pofgeneratingconst}}{\includegraphics[width=0.48\columnwidth, trim=2.4in .03in 2.7in 0in, clip=true]{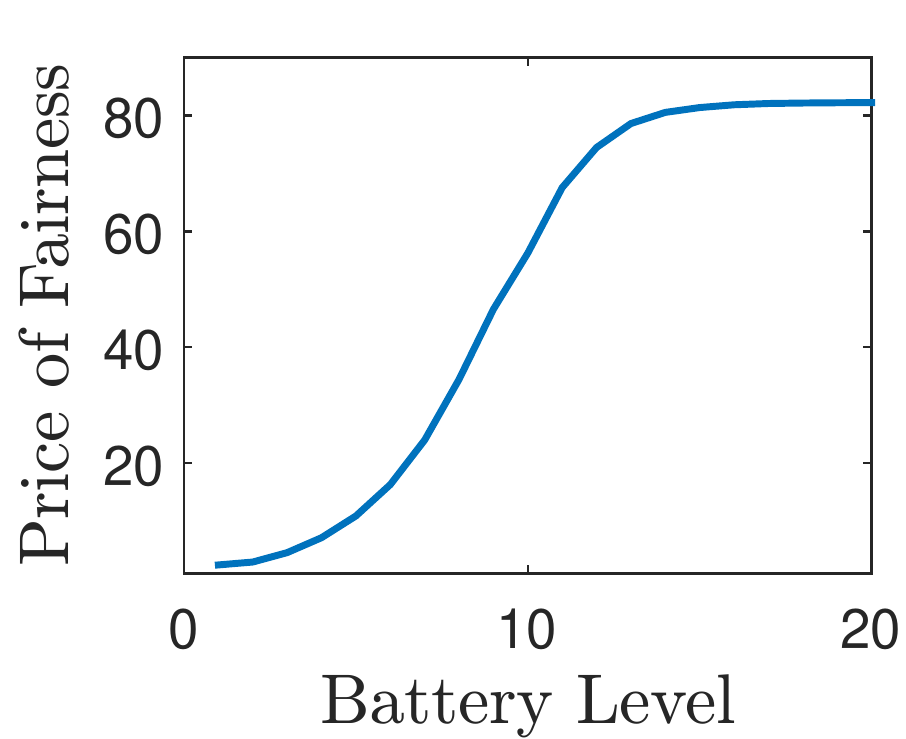}}
\captionsetup{font={small,sf},
    labelfont=bf,
    format=hang,    
    format=plain,
    margin=0pt,
    width=\columnwidth,}
\caption{These figures show the change in price of fairness as battery size increases when the sources are independent. $S = \{-1,1\}^n$ and $P_i$ denotes the transition matrix of source $i$.}
\end{figure}

An example of the phenomenon of unbounded PoF when all sources are net
generating is shown in~\cref{fig:pofgeneratingexp}. Note that the
above theorem does not show that the PoF for a given net generation
profile grows to infinity as $b_{\max} \ra \infty;$
see~\cref{fig:pofgeneratingconst}, where PoF seems to saturate with
increasing $b_{\max}.$

%% file: soft_fairness.tex
\section{Optimizing Fairness Subject to Efficiency}
\label{soft_fairness}

It is clear from Section~\ref{sec:pof} that there exist cases, in
which the Price of Fairness becomes unbounded as battery size
increases. In this section, we formulate a problem in which we try to
maximize fairness subject to efficiency, \ie, we pose being efficient
as a hard constraint and ask how fair can we get?

We first define the objective which captures the fairness in the
system. Our notion of fairness is to minimize the maximum over all
users, the amount a user draws minus the amount it supplies to the
battery. Recall the {\it net contribution} C$_i= \lim_{T \rightarrow
  \infty} \frac{1}{T}\sum_{t=0}^T \Ebb[A_i(t)]$ of user $i$ defined in
Section~\ref{sec:hard_fairness}. Now, by work conservation we have
under any policy, $\sum_{i \in [n]}{\rm C}_i = 0$. The system is said
to be \textit{completely fair} when ${\rm C}_i = 0 \ \forall i \in
[n]$. In general, this may not be the case.  Hence to maximize
fairness, we maximize the minimum of ${\rm C}_i$.

Here, by being \textit{efficient}, we mean being in the space of
policies which minimize the LLR without the hard fairness constraint.
In general, the optimal value of an objective, in this case fairness,
over an arbitrary set of policies is not easily computable. However,
an observation we made in Section~\ref{sec:pof} allows to formulate
this as a tractable problem.  Recall from Section~\ref{sec:pof} that
efficient policies are those that take actions in the space $A_e(x,b)$
when the state is $(x,b) \in\Sbf.$ As a consequence, the space of
efficient policies can be described by specifying only the space of
allowable actions. This allows us to formulate this problem as a CMDP.
%

Thus, we define the problem of maximizing fairness subject to
efficiency as follows
\begin{equation}\label{eqn:objfair}
\max_{\phi \in \Phi_e} \min_{i \in [n]} \lim_{T \rightarrow \infty} \frac{1}{T}\sum_{t=0}^T \Ebb[A_i(t)],
\end{equation}
where the maximization is over all $\phi$  that are efficient, \ie, that are constrained to take actions in $A_e(x,b)$ in every state $(x,b).$
This problem can also be formulated as a CMDP as follows    $$
    \maxproblemsmall{(F)}
    {\phi, \theta}
    {\theta}
                 {\begin{array}{r@{\ }c@{\ }l}
\theta  & \leq   & \lim_{T \rightarrow \infty} \frac{1}{T}\sum_{t=0}^T \Ebb[A_i(t)], \quad \forall i\in [n],  \\
\phi & \in & \Phi_e.
    \end{array}}
    $$

\subsection{Reduction to a Linear Program}\label{sec:lp}
The optimization problem (F) above is extremely complex since the
class of all history dependent policies is a highly unstructured
space. This makes it extremely hard to derive a policy or a battery
management algorithm. However, a remarkable simplification is
possible, which allows for the computation of such policies in an
efficient manner.

A Markov decision process is said to be unichain if, under any
stationary deterministic policy, the corresponding Markov chain
contains a single (aperiodic) ergodic class. This ensures the
existence of a unique stationary distribution, independent of the
initial distribution. Notice that under any policy $\phi \in \Phi_e$
the battery dynamics are identical and deterministic. Consequently,
under any policy $\phi \in \Phi_e$, the Markov chain $(X(t),b(t))$ has
a single ergodic class, making this CMDP unichain.  A well-known fact
about unichain CMDPs is that stationary randomized policies dominate
(see, \eg, Theorem 4.1 of \cite{altman99constrainedmdp}); thanks to
this, we can limit our search to stationary randomized
policies. Moreover, the CMDP admits an elegant equivalent linear
programming formulation (Theorem 4.3 of
\cite{altman99constrainedmdp}), where the optimization is not over
policies, but rather over \textit{occupation measures} (or probability
distributions) on the product of state and action spaces. We introduce
this formulation below.

An occupation measure $\rho$ is a probability distribution on $\Sbf \times \Abf_e$, where $\Abf_e = \bigcup_{(x,b) \in \Sbf}A_e(x,b)$. In the LP formulation, the occupation measure is required to satisfy,
\begin{align}\label{eqn:stat1}
\sum_{(x,b) \in \Sbf} \sum_{a\in A_e(x,b)} \rho(x,b,a) = 1, \quad \rho(x',b',a') & \geq 0,  \\
\hspace{-2mm}\sum_{(x,b)\in \Sbf} \sum_{a \in A_e(x,b)} \rho(x,b,a)( \delta_{(x',b')}(x,b) - P(x',b'|a,x,b))  & =0  \label{eqn:stat3} 
\end{align}
$\forall (x',b') \in \Sbf, a' \in \Abf_e$, where,
\begin{equation*}
 \delta_{(u,v)}(u',v') = \begin{cases}
                          1, &{\rm if} \ u = u' \ {\rm and} \ v = v',\\
                          0, &{\rm otherwise,}
                         \end{cases}
\end{equation*}
and $P(x',b'|a,x,b)$ is the probability of transition from state $(x,b)$ to state $(x',b')$ under action $a \in A_e(x',b')$. 
In our case, these dynamics are trivial:
\[P(x',b'|a,x,b) = P(x'|x)\I{b'=b+\sum_{i\in [n]} a_i },\]
and $P(x'|x)$ is the probability that the background process $X(t)$ transitions from $x$ to $x'$. \eqref{eqn:stat1} encodes that $\rho$ is a probability distribution on $\Sbf \times \Abf_e.$ \eqref{eqn:stat3} imposes that $\rho$ is an invariant distribution for the controlled Markov chain.

The equivalent LP formulation of (F) is then given by~\cite{altman99constrainedmdp}: 
$$    \maxproblemsmall{(LP)}
    {\rho, \theta}
    {\theta}
                 {\begin{array}{r@{\ }c@{\ }l}
\theta  & \leq   & \sum_{u = s^-_{i}}^{s^+_{i}}  \rho(a_i = u)u, \quad \forall i \in [n],  \\
\rho & \in & \rho_e.
    \end{array}}
    $$
Here, $\rho_e$ is the set of $\rho$ that satisfy \eqref{eqn:stat1}-\eqref{eqn:stat3}, 
and with a slight abuse of notation, we define 
\begin{align*}
\rho(a_i = u) &= \underset{(x,b,a) \in \Sbf \times \Abf_e: \ a_i = u}{\sum}\rho(x,b,a).
\end{align*}
The first constraint in (LP) arises from the fairness definition,
restated using occupation measures.  This formulation allows for
efficient computation of the optimal fairness. We next prove a result giving structural properties of efficient policies.

\begin{proposition}\label{prop:eff_prop}
Under any efficient policy, the battery dynamics remain the same. Moreover, the underlying marginal distribution over states ($\rho(x,b) \ \forall x,b \in \Sbf$) is the same across efficient policies.
\end{proposition}
\begin{IEEEproof}
Using E1, E2 and E3, we note that irrespective of what actions we choose for each agent, the sum total of all actions in $A_e(x,b)$ is the same for each pair $(x,b) \in \Sbf$. Since the battery dynamics are $b(t+1) = b(t) + \sum_{i}a_i(t)$, the battery dynamics stay the same across efficient policies.

Using the first claim, we can see that the transition matrix for $(x,b)$ is the same across all efficient policies. Thus, $\rho(x,b)$ is the same across all efficient policies.
\end{IEEEproof}

\newsavebox{\smlmatfairfir}
\savebox{\smlmatfairfir}{$\left[\begin{smallmatrix}0.5&0.5\\0.6&0.4\end{smallmatrix}\right]$}

\newsavebox{\smlmatfairsec}
\savebox{\smlmatfairsec}{$\left[\begin{smallmatrix}0.5&0.5\\0.8&0.2\end{smallmatrix}\right]$}

\newsavebox{\smlmatfairthir}
\savebox{\smlmatfairthir}{$\left[\begin{smallmatrix}0.7&0.3\\0.7&0.3\end{smallmatrix}\right]$}

\newsavebox{\smlmatfairfour}
\savebox{\smlmatfairfour}{$\left[\begin{smallmatrix}0.8&0.2\\0.4&0.6\end{smallmatrix}\right]$}

\newsavebox{\smlmatfairfiv}
\savebox{\smlmatfairfiv}{$\left[\begin{smallmatrix}0.6&0.4\\0.5&0.5\end{smallmatrix}\right]$}

\begin{figure}
\centering
\captionsetup{font={small,sf},
    labelfont=bf,
    format=hang,    
    format=plain,
    margin=0pt,
    width=\columnwidth,}
    \includegraphics[width = 0.8\columnwidth, trim=.25in .03in 0.75in 0.6in, clip=true]{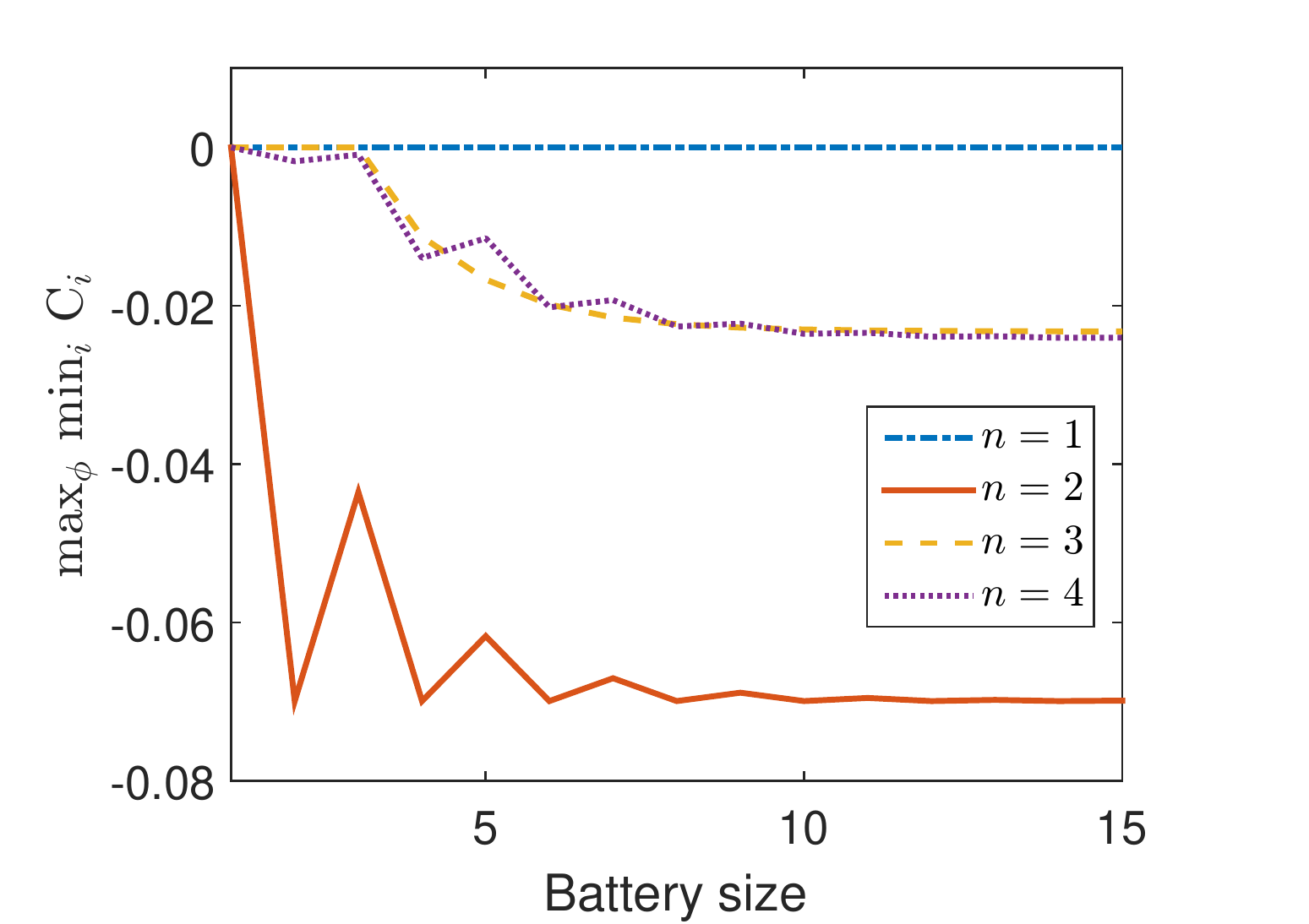}
 \caption{This figure shows how the maxmin fairness varies with increasing battery size. $S = \{-1,1\}^n$ and $P_i$ denotes the transition matrix of source $i$. We have the following sources and $n = k$ implies a system which has  first $k$ of these sources.\protect\footnotemark\\
  $P_1 = $ \usebox{\smlmatfairfir} , $P_2 = $ \usebox{\smlmatfairsec} , $P_3 = $ \usebox{\smlmatfairthir} , $P_4 = $ \usebox{\smlmatfairfour}.}
\label{fig:maxmin_fair_const}
 \end{figure}

\footnotetext{The irregularities for even number of users occur due to
  integer effects.}

We have computed the solution to (LP) and the plot of optimal
fairness with increasing battery size for different user
configurations is given in \cref{fig:maxmin_fair_const}. We have taken
the background processes of the users to be independent for these
evaluations. From this plot, we can observe that the fairness tends to
saturate on increasing battery size, \ie after a certain
threshold, the fairness value approaches a negative constant (except
for the case $n = 1,$ where the efficient policy is also fair). This
illustrates that increasing the battery size does not necessarily improve
the fairness of efficient policies. 

\cref{fig:maxmin_fair_const} indicates that even a relaxation of the
fairness constraint in~(P) to C$_i \geq -\delta$ for some $\delta > 0$
is not guaranteed to result in an efficient battery operation, even as
$b_{\max} \ra \infty.$ To illustrate this explicitly, in
\cref{fig:frontier3} and \cref{fig:frontier2}, we plot the fairness
efficiency frontier for 3 and 5 users, respectively. These plots show
how close we get to efficient operation by relaxing the fairness
constraint from C$_i \geq 0$ to C$_i \geq -\delta \ \forall i$. The
shaded region shows all those points which are achievable and the
frontier is the boundary for the achievable region. For the case when
an efficient policy is fair, the frontier would coincide with the
axes.

This reveals a fundamental conflict between efficiency and fairness: a
small relaxation of the fairness constraint does not necessarily
improve efficiency, and a small relaxation of the efficiency
requirement does not necessarily imply fairness for all users (this
follows from our previous observation of unbounded PoF).

\newsavebox{\smlmatfronfir}
\savebox{\smlmatfronfir}{$\left[\begin{smallmatrix}0.6&0.4\\0.6&0.4\end{smallmatrix}\right]$}

\newsavebox{\smlmatfronsec}
\savebox{\smlmatfronsec}{$\left[\begin{smallmatrix}0.9&0.1\\0.9&0.1\end{smallmatrix}\right]$}

\newsavebox{\smlmatfronthir}
\savebox{\smlmatfronthir}{$\left[\begin{smallmatrix}0.6&0.4\\0.5&0.5\end{smallmatrix}\right]$}

\newsavebox{\smlmatfronfour}
\savebox{\smlmatfronfour}{$\left[\begin{smallmatrix}0.7&0.3\\0.4&0.6\end{smallmatrix}\right]$}

\newsavebox{\smlmatfronfive}
\savebox{\smlmatfronfive}{$\left[\begin{smallmatrix}0.6&0.4\\0.9&0.1\end{smallmatrix}\right]$}

\newsavebox{\smlmatfronsix}
\savebox{\smlmatfronsix}{$\left[\begin{smallmatrix}0.6&0.4\\0.2&0.8\end{smallmatrix}\right]$}

\newsavebox{\smlmatfronsev}
\savebox{\smlmatfronsev}{$\left[\begin{smallmatrix}0.8&0.2\\0.8&0.2\end{smallmatrix}\right]$}

\newsavebox{\smlmatfroneight}
\savebox{\smlmatfroneight}{$\left[\begin{smallmatrix}0.51&0.49\\0.51&0.49\end{smallmatrix}\right]$}

\begin{figure}
\normalsize
\centering
\subcaptionbox{$ P_1  = $ \usebox{\smlmatfronfir}, $ P_2  =\ $ \usebox{\smlmatfronsec}, $ P_3  = $ \usebox{\smlmatfronthir} \label{fig:frontier3}}{\includegraphics[height=0.15\textwidth, trim=.15in 0.03in .25in .4in, clip=true]{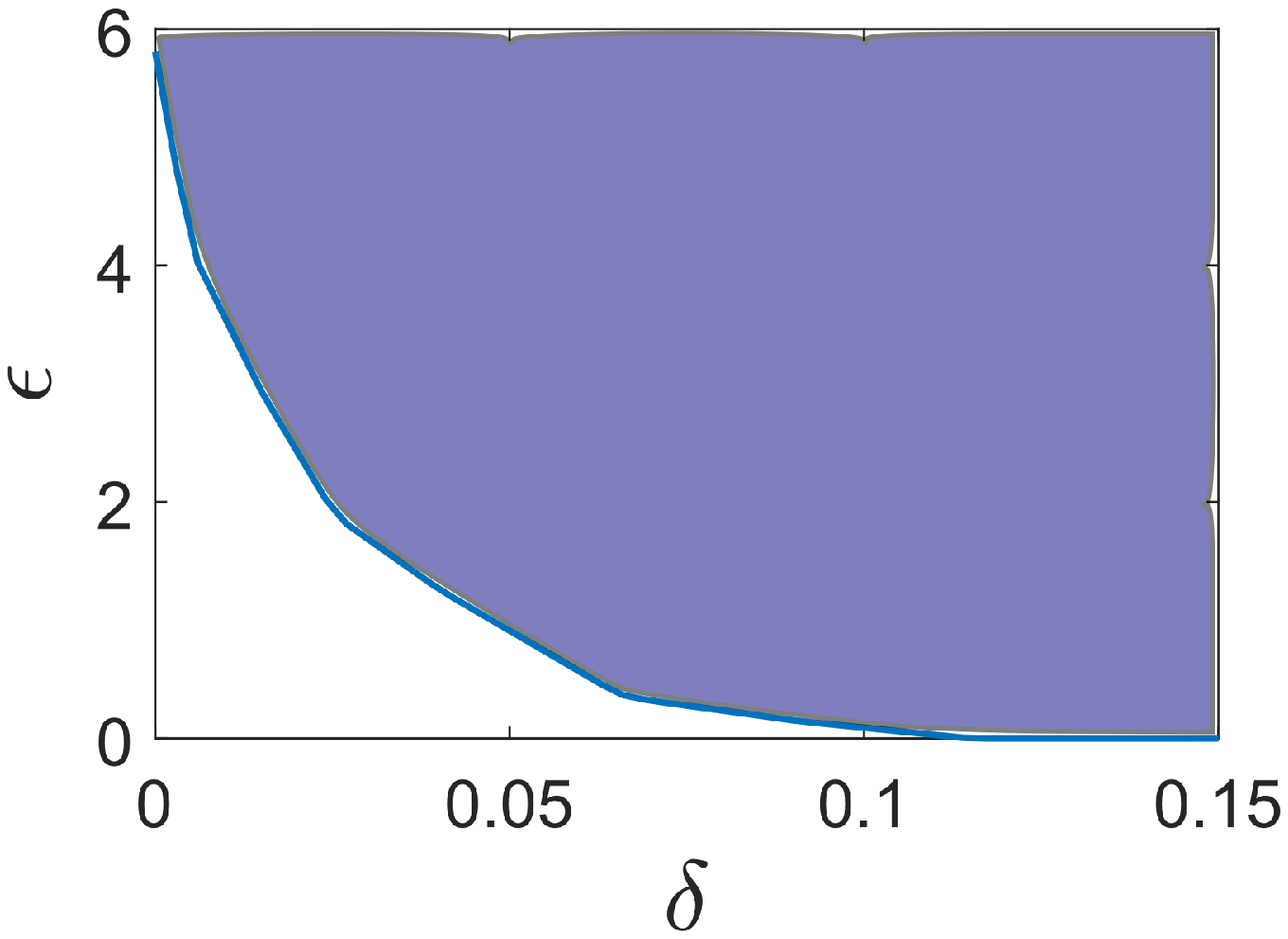}}\ \ %
\subcaptionbox{$ P_1  = $ \usebox{\smlmatfronfour}, $ P_2  = $ \usebox{\smlmatfronfive}, $ P_3  = $ \usebox{\smlmatfronsix}, $ P_4  = $ \usebox{\smlmatfronsev}, $ P_5  = $ \usebox{\smlmatfroneight}. \label{fig:frontier2}}{\includegraphics[height=0.15	\textwidth, trim=.1in .1in .25in .5in, clip=true]{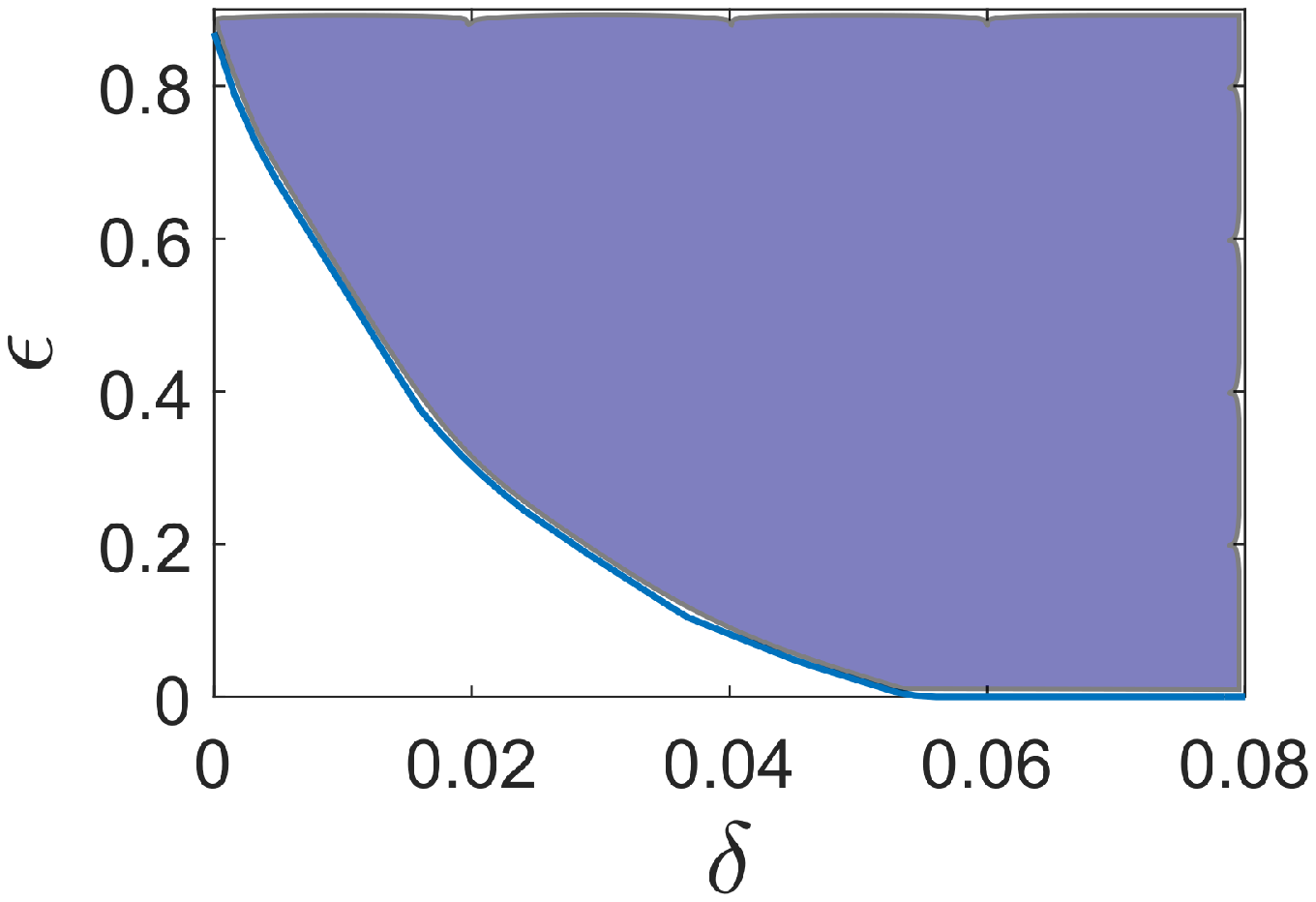}}
\captionsetup{font={small,sf},
    labelfont=bf,
    format=hang,    
    format=plain,
    margin=0pt,
    width=\columnwidth,}
\caption{These figures show the fairness efficiency frontier for 3 and 5 users respectively with $b_{\max} = 12$, $S = \{-1,1\}^n$ and $P_i$ denotes the transition matrix of source $i$. This optimal loss of load is denoted by $\LLR_\delta$. The plot is of $\epsilon$ v/s $\delta$, where $\epsilon = \frac{\LLR_\delta}{\LLR_e} - 1$. and the fairness constraints are $C_i \geq -\delta$.}
\end{figure}

To concretize our claim regarding the fundamental conflict between efficiency and fairness, we show in the following example that irrespective of the battery level, any efficient policy always has a non-zero maxmin fair value. We consider a setup of two users, wherein each user either generates unit power or demands unit power at each time instant. Let the transition probability matrix of  user $i$ be given by
\begin{equation}\label{eqn:transdef}
P_i = \left[\begin{matrix}p_i&1 - p_i\\1 - q_i & q_i\end{matrix}\right].
\end{equation}
Then, we define the generation probability of user $i$ by $\alpha_i = \frac{1 - q_i}{2 - p_i - q_i}$. Using this setup, we give the following proposition, which gives a class of examples where the fairness remains bounded away from zero under efficient policies. We denote the maxmin fair value under efficient policies by F$_e$.

\begin{proposition}
Consider two users with $S = \{-1,1\}^2$ and the transition probability matrix of user $i$ being $P_i$ (as in \cref{eqn:transdef}). If $\alpha_1 > \alpha_2$, we have F$_e = \alpha_2 - \alpha_1$ for all even battery sizes.
\end{proposition}
\begin{IEEEproof}
First, we note that $\rho(b = k) = 0$ for $k$ odd. This is because, for even battery size, we know that in the long run the battery spends only finite amount of time in the odd battery levels and spends an infinite amount of time in the even battery levels. To see this,  at each even battery  level, we either accept two units of power, accept one unit and give one unit, or supply two units of power. Thus the battery level doesn't change it's parity unless it hits a boundary. If the battery level is initially odd, we know that it will hit either the upper boundary ($b_{\max}$) or the lower boundary $(0)$ in finite time, both of which are even. After this, the battery level stays even.

Let $\alpha_i$ denote the generation probability of user $i$ and let $\alpha_i' = 1 - \alpha_i$. Let $\rho(\cdot,\cdot)$ denote the occupation measure of states of any efficient policy (Prop. \ref{prop:eff_prop}). We note that given that the policy is efficient, from E1,E2 and E3, the freedom in the action is only when the background state is $x = (-1,-1)$ and the battery level is $b = 1$ and when the background state is $x = (1,1)$ and the battery level is $b = b_{\max} - 1$. Let us assume that when the background state is $x = (-1,-1)$ and the battery level is $b = 1$, we choose action $a = (-1,0)$ with probability $s$ and action $a = (0,-1)$ with probability $1-s$. Also, assume that   when the background state is $x = (1,1)$ and the battery level is $b = b_{\max} - 1$, we choose action $a = (1,0)$ with probability $t$ and action $a = (0,1)$ with probability $1-t$. Thus, we have, 
\begin{align*}  
{\rm C}_1 - {\rm C}_2 = 
& 2 (\alpha_1 - \alpha_2)  + (2t - 1) \rho(x = (1,1),b = b_{\max} - 1) \\
& - (2s - 1) \rho(x = (-1,-1),b=1).
\end{align*}
Since $\rho(b = k) = 0$ for $k$ odd, we have 
${\rm C}_1 - {\rm C}_2 = 2 (\alpha_1 - \alpha_2).$
Thus, we have C$_1 = \alpha_1 - \alpha_2$ and C$_2 = \alpha_2 - \alpha_1$. For $\alpha_1 > \alpha_2$, C$_2 < $ C$_1$ and hence F$_e = \alpha_2 - \alpha_1$ for all even battery sizes.
\end{IEEEproof}

%% file: conclusion.tex
\section{Concluding Remarks and Future Work}
We studied the problem of scheduling of a shared energy storage among
multiple users by concentrating on the tradeoff between fairness of
use across users and the efficiency of use of the battery.  Our results
show that there are hard fundamental limits to fairness under
efficiency constraints and to efficiency under fairness
constraints. In particular, a larger battery size does not help in
achieving both fairness and efficiency simultaneously.  Our effort at
characterizing these tradeoffs was via optimization of efficiency
subject to fairness (and vice-versa). We have also simulated 
the \textit{fairness-efficiency frontier} of
levels of efficiency and fairness that cannot be jointly improved upon
under any policy or battery size. Our results also indicate that
sharing contracts and costing of batteries should be done with care since fair operations imply suboptimal battery utilization. Devising
an equitable sharing of costs and a market for transacting excess energy are fascinating directions of
future research.

%% file: appendix.tex
\section{Large battery asymptotics for a single user}\label{app:singlesource}

In this section, we consider the problem of minimizing the LLR
associated with a single user~$i$ operating a battery of size
$b_{\max}$ alone. We make the following remarks:
\begin{enumerate}
\item The net generation process $X_i(\cdot)$ corresponding to
  user~$i$ is a functional of the DTMC $X(\cdot).$
\item An elementary energy conservation argument shows that any policy
  of battery operation for user~$i$ is fair with respect to the
  fairness notion~\eqref{eqn:fair-cons}.
\item The LLR of user~$i$ is minimized by a simple greedy policy:
  Suppose that the net generation is $x_i$ and the battery occupancy
  equals $b.$ If $0 \leq b+x \leq b_{\max}$, set action $a = x.$ If $b
  + x < 0$, set $a = -b,$ and if $b + x > b_{\max}$, set $a = b_{\max}
  - b$. Under this policy, the battery evolution is given by 
  \begin{equation}
    \label{eq:single_battery_evol}
    B(t+1)   = [B(t) + X_i(t)]_{[0,b_{\max}]}, 
  \end{equation}
  where $[z]_{[0,b_{\max}]} = \min(\max(z,0),b_{\max}).$
\end{enumerate}

The main result of this section is that if $\Delta_i>0,$ then
$\LLR_{o,i},$ the optimal LLR for user~$i,$ decays exponentially
with the battery size~$b_{\max}.$

\begin{theorem}
\label{lem:netgensingle}
For user~$i$ operating a battery of size $b_{\max}$ in a standalone
fashion, if $\Delta_i>0,$ then  
\[
\lim_{b_{\max} \rightarrow \infty}\frac{\log(\LLR_{o,i})}{b_{\max}} = -\lambda_i,
\]
where $\lambda_i \in (0,\infty).$
\end{theorem}

Theorem~\ref{lem:netgensingle} is a consequence of the following
lemmas.

\begin{lemma}
  \label{lemma:llrsingle}
  Let the setting of Theorem~\ref{lem:netgensingle} hold. Under the
  battery evolution given by \eqref{eq:single_battery_evol},
  $$\lim_{b_{\max} \rightarrow \infty}\frac{\log(\LLR_{o,i})}{b_{\max}} 
  = \lim_{b_{\max} \rightarrow \infty}\frac{\log \prob{B=0}}{b_{\max}}.$$
where $B$ is distributed according to the steady state distribution of $B(\cdot)$
\end{lemma}

\begin{IEEEproof}
  It is easy to see that $\LLR_{o,i} \leq \prob{B = 0}.$ It therefore
  suffices to show that $\LLR_{o,i} \geq c \prob{B = 0}$ for some $c
  \in (0,1).$ This is a direct consequence of our assumption that
  $\prob{X(t+1) = s\ |\ X(t) = s} > 0$ for all $s \in S,$ which
  ensures that each time the battery drains to zero, there is positive
  probability of a loss of load. This argument can be easily
  formalized using the renewal reward theorem.
\end{IEEEproof}

\begin{lemma}
\label{lemma:battery_zero_asymptotics}
Let the setting of Theorem~\ref{lem:netgensingle} hold. Under the
battery evolution given by
\eqref{eq:single_battery_evol}, 
\begin{equation*}
  \label{eq:batt_zero_asymptotics}
  \lim_{b_{\max} \rightarrow
    \infty}\frac{\log \prob{B=0}}{b_{\max}} = -\lambda_i,
\end{equation*}
where $\lambda_i \in (0,\infty).$
\end{lemma}
\begin{IEEEproof}
  We analyse the asymptotics of $\prob{b = 0}$ using the
  \emph{reversed} system \cite{Mitra88}, which is obtained by
  interchanging the role of generation and demand. In other words,
  $X_i^{\rev}(t) = -X_i(t)$ (we use the superscript $\rev$ to
  represent quantities in the reversed system). Thus, $\Delta_i^{\rev}
  = -\Delta_i < 0.$ It is not hard to see that $$\prob{B = 0} =
  \prob{B^{\rev} = b_{\max}}.$$

  The battery evolution in the reversed system is equivalent to the
  buffer evolution in a finite buffer Markov modulated queue with
  negative drift, for which logarithmic asymptotics of the form $$
  \lim_{b_{\max} \rightarrow \infty}\frac{\log
    \prob{B^{\rev}=b_{\max}}}{b_{\max}} = -\lambda$$ are known; see
  \cite[Section~6.5]{BigQueues} and \cite{Toomey98}.
\end{IEEEproof}

%% file: root-new.bbl
\begin{thebibliography}{10}
\providecommand{\url}[1]{#1}
\csname url@samestyle\endcsname
\providecommand{\newblock}{\relax}
\providecommand{\bibinfo}[2]{#2}
\providecommand{\BIBentrySTDinterwordspacing}{\spaceskip=0pt\relax}
\providecommand{\BIBentryALTinterwordstretchfactor}{4}
\providecommand{\BIBentryALTinterwordspacing}{\spaceskip=\fontdimen2\font plus
\BIBentryALTinterwordstretchfactor\fontdimen3\font minus
  \fontdimen4\font\relax}
\providecommand{\BIBforeignlanguage}[2]{{%
\expandafter\ifx\csname l@#1\endcsname\relax
\typeout{** WARNING: IEEEtran.bst: No hyphenation pattern has been}%
\typeout{** loaded for the language `#1'. Using the pattern for}%
\typeout{** the default language instead.}%
\else
\language=\csname l@#1\endcsname
\fi
#2}}
\providecommand{\BIBdecl}{\relax}
\BIBdecl

\bibitem{chadha2019fair}
K.~N. Chadha, A.~A. Kulkarni, and J.~Nair, ``Efficiency fairness tradeoff in
  battery sharing,'' in \emph{Decision and Control (CDC), 2019 IEEE 58nd Annual
  Conference on}.\hskip 1em plus 0.5em minus 0.4em\relax IEEE, 2019, p. pages.

\bibitem{kalathil2019sharingelec}
\BIBentryALTinterwordspacing
D.~Kalathil, C.~Wu, K.~Poolla, and P.~Varaiya, ``The sharing economy for the
  electricity storage,'' \emph{{IEEE} Transactions on Smart Grid}, vol.~10,
  no.~1, pp. 556--567, jan 2019. [Online]. Available:
  \url{https://doi.org/10.1109/tsg.2017.2748519}
\BIBentrySTDinterwordspacing

\bibitem{chakraborty2018coalationalsharing}
\BIBentryALTinterwordspacing
P.~Chakraborty, E.~Baeyens, K.~Poolla, P.~P. Khargonekar, and P.~Varaiya,
  ``Sharing storage in a smart grid: A coalitional game approach,''
  \emph{{IEEE} Transactions on Smart Grid}, pp. 1--1, 2018. [Online].
  Available: \url{https://doi.org/10.1109/tsg.2018.2858206}
\BIBentrySTDinterwordspacing

\bibitem{wu2016communitystorage}
\BIBentryALTinterwordspacing
C.~Wu, J.~Porter, and K.~Poolla, ``Community storage for firming,'' in
  \emph{2016 {IEEE} International Conference on Smart Grid Communications
  ({SmartGridComm})}.\hskip 1em plus 0.5em minus 0.4em\relax {IEEE}, nov 2016.
  [Online]. Available: \url{https://doi.org/10.1109/smartgridcomm.2016.7778822}
\BIBentrySTDinterwordspacing

\bibitem{zhou2018windstoragemdp}
\BIBentryALTinterwordspacing
Y.~H. Zhou, A.~Scheller-Wolf, N.~Secomandi, and S.~Smith, ``Managing wind-based
  electricity generation in~the~presence of storage and transmission
  capacity,'' \emph{Production and Operations Management}, nov 2018. [Online].
  Available: \url{https://doi.org/10.1111/poms.12946}
\BIBentrySTDinterwordspacing

\bibitem{kim2011optenergy}
\BIBentryALTinterwordspacing
J.~H. Kim and W.~B. Powell, ``Optimal energy commitments with storage and
  intermittent supply,'' \emph{Operations Research}, vol.~59, no.~6, pp.
  1347--1360, dec 2011. [Online]. Available:
  \url{https://doi.org/10.1287/opre.1110.0971}
\BIBentrySTDinterwordspacing

\bibitem{korpaas2003opersizingstorage}
\BIBentryALTinterwordspacing
M.~Korpaas, A.~T. Holen, and R.~Hildrum, ``Operation and sizing of energy
  storage for wind power plants in a market system,'' \emph{International
  Journal of Electrical Power {\&} Energy Systems}, vol.~25, no.~8, pp.
  599--606, oct 2003. [Online]. Available:
  \url{https://doi.org/10.1016/s0142-0615(03)00016-4}
\BIBentrySTDinterwordspacing

\bibitem{vandevan2012optcontenergy}
P.~M.~van~de Ven, N.~Hegde, L.~Massoulié, and T.~Salonidis, ``Optimal control
  of end-user energy storage,'' \emph{IEEE Transactions on Smart Grid}, vol.~4,
  03 2012.

\bibitem{altman99constrainedmdp}
E.~Altman, \emph{Constrained Markov Decision Processes}, 1999.

\bibitem{PutermanMDP}
M.~L. Puterman, \emph{Markov decision processes: discrete stochastic dynamic
  programming}.\hskip 1em plus 0.5em minus 0.4em\relax John Wiley \& Sons,
  2014.

\bibitem{Mitra88}
D.~Mitra, ``Stochastic theory of a fluid model of producers and consumers
  coupled by a buffer,'' \emph{Advances in Applied Probability}, vol.~20,
  no.~3, pp. 646--676, 1988.

\bibitem{BigQueues}
A.~J. Ganesh, N.~O'Connell, and D.~J. Wischik, \emph{Big queues}.\hskip 1em
  plus 0.5em minus 0.4em\relax Springer, 2004.

\bibitem{Toomey98}
F.~Toomey, ``Bursty traffic and finite capacity queues,'' \emph{Annals of
  Operations Research}, vol.~79, pp. 45--62, 1998.

\end{thebibliography}
